\documentclass[11pt]{article}

\usepackage{amssymb,amsmath,a4wide}
\usepackage{textcomp}
\usepackage{enumitem}
\usepackage[colorlinks=true,urlcolor=blue,citecolor=blue,linkcolor=blue]{hyperref}

\usepackage[T1]{fontenc}

\newcommand{\R}{\mathbb{R}}
\newcommand{\LL}{\mathbb{L}^3}
\renewcommand{\L}{\mathbb{L}}
\newcommand{\I}{\mathbf{I}}

\newtheorem{definition}{Definition}[section]
\newenvironment{pf}
{\par\noindent\textit{Proof:}\,}{$\Box$\par}

\newtheorem{thm}{Theorem}[section]

\newtheorem{rmk}{Remark}[section]
\newtheorem{corollary}{Corollary}[section]

\newtheorem{example}{Example}[section]

\begin{document}

\title{A Zoo of Translating Solitons on a Parallel Light-like Direction in Minkowski 3-Space}

\author{Erdem Kocaku\c{s}akl\i \\ 
\small Department of Mathematics, Faculty of Science, University \\
\small of Ankara Tandogan, Ankara, TURKEY\\
kocakusakli@ankara.edu.tr \and Miguel Ortega \\
\small Institute of Mathematics, Department of Geometry and Topology,\\
\small University of Granada, Granada, SPAIN\\
miortega@ugr.es}

\date{}%\today}
% \address[2]{Institute of Mathematics, Department of Geometry and Topology, University of Granada, 18071 Granada, Spain}
%\ead{miortega@ugr.es}
\maketitle

\begin{abstract} We deal with solitons of the mean curvature flow. The definition of \textit{translating solitons on a light-like direction} in Minkowski 3-space is introduced.  Firstly, we classify those which are graphical, \textit{translation surfaces}, obtaining space-like and time-like, entire and not entire, complete and incomplete examples. Among them, all our  time-like examples are incomplete. The second family consists of those  which are invariant by a 1-dimensional subgroup of parabolic motions, i.e, with light-like axis. The classification result implies that all examples of this second family have singularities. 
\end{abstract}



\section{Introduction}

Hypersurfaces in Euclidean space which evolve along the mean curvature flow have been widely studied. Of particular interest are those called \textit{translating solitons}, which are those whose mean curvature vector satisfies the following equation 
\begin{equation*}
\vec{H}=\vec{K}^{\perp}
\end{equation*}
where $\perp $ and $\vec{K}$ denote projection on  the normal bundle and a unit vector field, respectively. A much more general, but very weak, definition can be found in \cite{ALR}, where virtually no restriction on $\vec{K}$ is set (see also \cite{LiraM}). Among many papers, we can select \cite{ALR}, \cite{LiraM}, \cite{HIMW},  which study translating solitons in Riemannian manifolds. However, the theory seems to be less developed in Lorentzian geometry, although relevant results can be found in \cite{ding} and \cite{SX}. According to them, translating solitons in Minkowski space are mainly studied when the vector $\vec{K}$ is parallel and time-like. 

In this paper, we wish to introduce a new family, namely  translating solitons of the mean curvature flow such that $\vec{K}$ is a \textit{light-like parallel vector}. More precisely. 
\begin{definition}\label{Tlight-like} Given a parallel  light-like vector $\vec{K}$ in Minkowski 3-space $\LL$, a non-degenerate immersion $\psi:M\to \LL$ will be called a {\em translating soliton on the light-like direction} 
$\vec{K}$ if its mean curvature vector $\vec{H}$ satisfies $\vec{H}=\vec{K}^{\perp}$, where $\perp$ means the orthogonal projection on the normal bundle.
\end{definition}
\noindent This definition makes sense because the induced metric is not degenerate, which implies that $\vec{H}$ will not be light-like at any point. 

We write the standard flat metric in Minkowski 3-space in a suitable way for our needs, namely  $\langle,\rangle=-2dxdy+dz^{2}$. That is to say, we are considering a basis $\mathcal{B}=\{\vec{x},\vec{y},\vec{z}\}$ such that $\vec{x}$, $\vec{y}$ are light-like, future pointing, and satisfying the normalizing condition $\langle\vec{x},\vec{y}\rangle =-1$.  It is important to recall that any two parallel  light-like vectors are linked by an isometry of $\L^3$, because the light cone is invariant by rotations, boosts and a few reflections. This means that we can reduce to the case $\vec{K}=\vec{x}$. We will study two families. 

Firstly, Section \ref{graphical} is devoted to studying graphical surfaces, that is to say, the ones that admit a parametrization  $\psi:\Omega\subset\R^2\to \L^3$, $\psi(y,z)=(u(y,z),y,z)$, where $u:\Omega\to\R$. In Theorem \ref{graphical}, we will classify those which, in addition, are \textit{translation surfaces}, i.~e., for some smooth functions $a$ and $b$, then $u(y,z)=a(y)+b(z)$.  Translation surfaces in Euclidean space were introduced by S. Lie (see \cite{darboux}, also \cite{LM}). Needless to say, the same definition can be easily set in Minkowski space. We study space-like and time-like surfaces, obtain four types, which are flat by chance. Note that \cite{ding} and \cite{SX} paid attention to entire examples.  We later study the completeness of the four cases, showing entire and not entire, complete and incomplete surfaces, along Corollaries \ref{tipo-I}, \ref{tipo-II-completeness}, \ref{tipo-III} and \ref{type1}. We should remark that the standard techniques to show the completeness of space-like surfaces in Minkowski space do not work in our setting (see Remark \ref{novan}).

Secondly, there are 1-dimensional subgroups of parabolic isometries of $\L^3$, whose axis is light-like. We will reduce to the well-adapted case when the \textit{rotation axis} is spanned by the vector $\vec{x}$ (see  Section \ref{eje} for more details.) That is to say, we study surfaces obtained by letting this subgroup of isometries act on a suitable profile curve.  We classify in Theorem \ref{parabolic} those translating solitons which are invariant by this subgroup of isometries. 

\section{Preliminaries}

Given a smooth manifold $M$, assume a family of smooth immersions in a semi-Riemannian manifold $(M,g)$, $F_{t}:M\rightarrow M\times $ $\mathbb{R}$, $t\in \left[ 0,\delta \right)$, $\delta >0$, with mean curvature vector $\vec{H}_{t}$. The initial immersion $F_{0}$ is called a \textit{solution to the Mean Curvature Flow} (MCF), up to local diffeomorphism, if the following equation holds
\begin{equation*}
\left( \frac{d}{dt}F_{t}\right) ^{\perp}=\vec{H}_{t},
\end{equation*}
where $\perp$ means the orthogonal projection on the normal bundle. If an immersion $F:M\rightarrow \LL$ satisfies the condition $\vec{H}=\vec{K}^{\perp}$, then it is possible to define the forever flow 
$\psi :M\times \R \rightarrow \LL$, $\psi \left( p,t\right) =F\left( p\right) +t\vec{K}$, and clearly, 
\begin{equation*}
\left( \frac{d}{dt}F_{t}\right) ^{\perp}=\vec{K}^{\perp}= 
\vec{H}.
\end{equation*}

We recall the basic theory of surfaces in Minkowski 3-space. See \cite{LO} for details. Let $\psi:M\to \LL$ be an immersion of a surface $M$ in Minkowski 3-space. We assume that the induced metric $\I=\psi^*\langle,\rangle$ is not degenerate, that is to say, it is either Riemannian or Lorentzian. This metric is usually known as the \textit{first fundamental form}.  Let $N$ be a unit normal vector field on $M$. Therefore, $\varepsilon=\langle N,N\rangle=\pm 1$ is a constant function on $M$. Given $\nabla^o$ and $\nabla$ the Levi-Civita connection of $\LL$ and the induced connection on $M$, respectively, we know
\[ \nabla^o_XY = \nabla_XY +\varepsilon \langle AX,Y\rangle N, \quad X,Y\in TM,
\]
where $A$ is the shape operator of $N$. The \textit{second fundamental form} $\sigma$ is defined as 
\[ \sigma(X,Y)= \varepsilon \langle AX,Y\rangle N, \quad X,Y\in TM.
\]
Our definition of the \textit{mean curvature vector} is 
\[ \vec{H}= \varepsilon H N=\mathrm{trace}_{\I}(\sigma), \quad H=\mathrm{trace}_{\I}(A).
\]
The function $H:M\to \R$ is called the \textit{mean curvature} (function) of $M$. As in \cite{LO}, if  $X$ is a local parametrization $X:U\subset \R^2\to \LL$, $X=X(u,v)$, then  $B=(X_u,X_v)$ is a local basis of the tangent plane at each point of $X(U)$. The \textit{coefficients} of the first fundamental form are
\[ E=\langle X_u,X_u\rangle,\ F=\langle X_u,X_v\rangle,\ G=\langle X_v,X_v\rangle, 
\]
so that the matricial expression is $\left(\begin{smallmatrix} E & F \\ F & G \end{smallmatrix}\right)$.  Let $\left(\begin{smallmatrix} e & f \\ f& g\end{smallmatrix}\right)$ be the matricial expression of $\sigma$ with respect to $B$. Since $\langle AX,Y\rangle=\langle \nabla_X^oY,N\rangle$, then 
\begin{gather*}
e  = \langle AX_u,X_u\rangle = \langle N,X_{uu}\rangle, \
f  = \langle AX_u,X_v\rangle = \langle N,X_{uv}\rangle, \ 
g  = \langle AX_v,X_v\rangle = \langle N,X_{vv}\rangle. 
\end{gather*}
The shape operator can be computed by 
\begin{equation}\label{shape}
A\equiv \begin{pmatrix} E & F \\ F & G \end{pmatrix}^{-1}
\begin{pmatrix} e & f \\ f & g \end{pmatrix}.
\end{equation}
With this, the expressions of the mean and Gaussian curvatures are
\begin{equation}\label{formulaH}
H = \frac{Eg-2Ff+Ge}{EG-F^2}, \quad 
K=\dfrac{eg-f^2}{EG-F^2}. 
\end{equation}

\section{Graphical Translating Solitons on a Light-like Direction}\label{graphical}

Given a domain $\Omega\subset\R^2$, let us take a parametrization of a non-degenerate surface  
\begin{equation*}
\psi :\Omega \rightarrow \mathbb{L}^{3},\ 
\psi(y,z) =\left( u(y,z),y,z\right).
\end{equation*}
The partial derivatives of $\psi (y,z)$ are $\psi_y=(u_y,1,0)$ and $\psi_z=(u_z,0,1).$ 
The coefficients of the first fundamental form $\I=\psi^*\langle,\rangle$ are 
\begin{equation*}
\left(\begin{matrix} E & F \\ F & G \end{matrix}\right) =
\left(\begin{matrix} -2u_y & -u_z \\ 
-u_z & 1 \end{matrix}\right).
\end{equation*}
In addition, it is Riemannian if $EG-F^2=-2u_y-u_z^{2}>0$. Since we are assuming that the surface is not degenerate, the following function is constant,  $\varepsilon=\mathrm{sign}(2u_y+u_z^{2})=\pm 1$. A unit normal vector field is given by 
\[ N=\frac{1}{W}(-u_y,1,u_z), \ \textrm{where} \ W=\sqrt{\varepsilon\big(2u_y+u_z^2\big)}>0. 
\]
Needless to say, $\langle N,N\rangle=\varepsilon$. 
The coefficients of the second fundamental form are
\[ e=\langle N,\psi_{yy}\rangle = \frac{-u_{yy}}{W},\quad 
f =\langle N,\psi_{yz}\rangle= \frac{-u_{yz}}{W}, \quad 
g=\langle N,\psi_{zz}\rangle=\frac{-u_{zz}}{W}. 
\]
We  compute the mean curvature function
\[
H =  \frac{-2u_y\,u_{zz} + 2u_z\,u_{yz} + u_{yy} }{W(2u_y+u_z^2)}.
\]
On the other hand, by using the definition of the translating soliton on the
light-like direction $\vec{K}=\vec{x}$, we obtain 
$\varepsilon H N = \vec{H} = \vec{x}^{\perp}$, so that $H=g(\vec{H},N)=g(\vec{x}^{\perp},N)
=g(\vec{x},N)=-1/W$. We finally obtain the following PDE, which characterizes our graphical translating solitons in the light-like direction $\vec{K}=\vec{x}$, 
\begin{equation}
u_{yy}+2u_zu_{yz}-2u_yu_{zz}+2u_y+u_z^{2}=0.  \label{PDE}
\end{equation}

Now, we consider $\psi$ as a \textit{translation surface}, which means  $u(y,z) =a(y) +b(z)$, where $a$ and $b$ are smooth functions. This is equivalent to $u_{yz}=0$ everywhere. 
So, our partial differential equation turns into the following,
\begin{equation}
a''(y)-2a'(y)b''(z)+2a'(y)+\left( b'(z) \right) ^{2}=0.  \label{mainpde}
\end{equation}

\noindent \textsc{Case 1.} $a''(y)=0$ in an interval. Then  $a(y)=a_1y+a_0$ for some $a_0,a_1\in\R$. In \eqref{mainpde}, we obtain 
\[-2a_1b''(z)+2a_1+(b'(z))^2=0.\] 
If $a_1=0$, then $b'(z)=0$ in an interval, that is to say, $u(y,z)=a_0+b_1$. But we are discarding this case because we need $2u_y+u_z^2\neq 0$. We put $\varphi =b'$, so that 
\[ \dfrac{\varphi'(z)}{(\varphi(z))^2+2a_1}=\dfrac{1}{2a_1}.
\]
\noindent\textsc{Case 1.1. $a_1=2\lambda^2>0$:} 
By integrating, and recalling that $b'(z)=\varphi (z)$, we get for some $z_0\in\R$, 
\begin{equation*}
b(z)=-4\lambda ^{2}\log_e \left\vert \cos \left( \frac{z-z_0}{2\lambda }
\right) \right\vert +b_{0},
\end{equation*}
where $b_{0}$ is an integration constant. We reach to the following solution 
\begin{equation}
u(y,z)=2\lambda ^{2}y-4\lambda ^{2}\log_e \left\vert \cos \left( \frac{z-z_0}{
2\lambda }\right) \right\vert +b_{0},\quad \lambda >0,\, b_{0},z_0\in \mathbb{R}.  \label{solution3}
\end{equation}

\noindent\textsc{Case 1.2. $a_1=-2\lambda^2<0$:} If we integrate both sides, then we obtain the following,
\[\log_e \left\vert \dfrac{\varphi(z)-2\lambda}{\varphi(z)+2\lambda}\right\vert 
=\dfrac{z-z_0}{\lambda}.
\]
From here, we have to discuss two cases, namely  positive and negative:
\begin{gather*}
b'(z)=\varphi(z)=2\lambda \coth\left(\dfrac{z-z_0}{2\lambda}\right),\ 
b'(z)=\varphi(z)=2\lambda \tanh\left(\dfrac{z-z_0}{2\lambda}\right). 
\end{gather*}
Our solutions are
\begin{align}
u(y,z)=& -2\lambda ^{2}y-4\lambda ^{2}\log_e \left\vert \sinh \left( \frac{
z-z_{0}}{2\lambda }\right) \right\vert +a_{0},  \\
u(y,z)=& -2\lambda ^{2}y-4\lambda ^{2}\log_e \left( \cosh \left( \frac{z-z_{0}}{
2\lambda }\right) \right) +a_{0}.
\end{align}

\noindent \textsc{Case 2.} $a''(y)\neq 0$ in an interval.\\
\noindent \textsc{Case 2.1.} We differentiate equation \eqref{mainpde} with respect to $y$: 
\[b''(z) =\frac{a'''(y) +2a''(y)}{2a''(y) }=b_2\in \mathbb{R}.
\]
Therefore, $b(z)=b_2z^{2}/2+b_{1}z+b_{0}$ for some integration constants $b_1,b_0\in\R$. 
We return to \eqref{mainpde}, 
\begin{equation*}
a''(y)-2a'(y)b_2+2a'(y)+\left( b_2z+b_1 \right) ^{2}=0.
\end{equation*}
This must hold on some intervals. This readily implies that $b_2=0$. Next, we
obtain 
\begin{equation*}
a''(y)+2a'(y)+b_{1}^{2}=0.
\end{equation*}
(Note that this implies $a'''(y)+2a''(y)=0$). The general solution to this equation is 
\begin{equation*}
a(y)=a_{1}e^{-2y}-\frac{b_{1}^{2}}{2}y+a_{0},\quad a_{0},a_{1}\in \mathbb{R}.
\end{equation*}
Therefore, a  family of solutions to \eqref{PDE} is 
\begin{equation}
u:\mathbb{R}^{2}\rightarrow \mathbb{R},\quad 
u(y,z) = a_{1}e^{-2y}-\frac{b_{1}^{2}}{2}y+b_{1}z+b_{0},\quad a_{1},b_{1},b_{0}\in \mathbb{R}.
\label{solution1}
\end{equation}

\noindent\textsc{Case 2.2.} Similarly to the previous idea, we differentiate
equation \eqref{mainpde} with respect to $z$: 
\[ b'(z)b''(z) =a'(y)b'''(z).\]

\noindent \textsc{Case 2.2.1.} If $b'''(z)=0$ in an interval, then  $b(z)=b_2z^2+b_{1}z+b_{0}$, $b_{0},b_{1},b_2\in \mathbb{R}$. This is the same function as in the previous Case 2.1, so that we obtain again solution \eqref{solution1}.\\ 

\noindent \textsc{Case 2.2.2.} If $b'''(z)\neq 0$ in an interval, then $\dfrac{b'(z)b''(z)}{b^{\prime
\prime \prime }(z)}=a'(y)=a_1$, where $a_1\in \mathbb{R}$ is constant. If $a_1=0$, $b'(z)b''(z)=0$, and we are again in Case 2.2.1. So, we assume $a_1\neq 0$. Obviously, 
\begin{equation*}
a(y)=a_1y+a_{0},\quad a_1,a_{0}\in \mathbb{R}.
\end{equation*}
Moreover, we have the following equality $b'(z)b''(z)=a_1b'''(z)$, which implies 
\begin{equation*}
\frac{\left( b'(z)\right) ^{2}}{2}=a_1\,b''(z)+b_0.
\end{equation*}
We return to \eqref{mainpde}, so that $0=2a_1+2b_0$, that is to say, $b_0=-a_1$. This was already discussed in Case 1.  \\

\begin{thm}
\label{convolution-surfaces} \label{psi} Consider $\LL$ with the usual metric written as $g=-2dxdy+dz^2$. Let $\psi:\Omega\subset\R^2\to\LL$, $\psi (y,z) =\left(
u(y,z),y,z\right)$ be a translation surface, i.e., we write $u(y,z) =a(y)+b(z)$,
where $a$ and $b$ are smooth functions. If $\psi$ is a translating soliton on
the parallel light-like direction $\vec{K}=\vec{x}=\left(1,0,0\right)$, then $u$ is one of the following:
\begin{enumerate}[label={\Roman*}.]
\item \label{caso1} $u(y,z)=-2\lambda ^{2}y+4\lambda ^{2}\log_e \left( \cosh
\left( \dfrac{z-z_0}{2\lambda }\right) \right) +a_{0}$, $\lambda >0$, $a_{0}$, $z_0\in \mathbb{R}$.

\item \label{caso2} $u(y,z)=a_{1}e^{-2y}-\dfrac{b_{1}^{2}}{2}y+b_{1}z+b_{0}$, 
$a_{1},b_{1},b_{0}\in \mathbb{R}$, $a_{1}\neq 0$;

\item \label{caso3} $u(y,z)=2\lambda ^{2}y-4\lambda ^{2}\log_e \left\vert \cos
\left( \dfrac{z-z_0}{2\lambda }\right) \right\vert +b_{0}$, $\lambda >0$, $b_{0}$, $z_0\in \mathbb{R}$;

\item \label{caso4} $u(y,z)=-2\lambda ^{2}y+4\lambda ^{2}\log_e \left\vert
\sinh \left( \dfrac{z-z_0}{2\lambda }\right) \right\vert +a_{0}$, $\lambda
>0$, $a_{0}$, $z_0\in \mathbb{R}$;
\end{enumerate}
\end{thm}
We will call them \textit{of type I, II, III and IV}, respectively. A straightforward computation shows that $K=0$ in all cases.
\begin{corollary} \label{flat}  All examples in Theorem \ref{convolution-surfaces} are flat.
\end{corollary}

Next, we want to study the completeness of these surfaces. 

\begin{corollary}\label{tipo-I}
\label{type3-4} Any inextensible solution of type $I$ %\ref{caso1} 
in Theorem \ref{convolution-surfaces} is space-like, entire and complete.  
\end{corollary}
\begin{pf} Given $a_{0},z_{0}\in \mathbb{R}$, $\lambda>0$, we consider 
\[u(y,z)=-2\lambda ^{2}y+4\lambda ^{2}\log_e\left( \cosh \left( \dfrac{z-z_{0}}{2\lambda }\right) \right) +a_{0}.\] Needless to say, it can be extended to $u:\mathbb{R}^2\to\R$. We compute 
\begin{gather*}
u_y=-2\lambda ^{2},\quad  u_z=2\lambda \tanh \left( \frac{z-z_{0}}{2\lambda }\right),\\
E=-2u_y=4\lambda^2, \ F=-u_z=-2\lambda \tanh \left( \frac{z-z_{0}}{2\lambda }\right), \ G=1,\\
EG-F^2=\frac{4\lambda^2}{\cosh^2\left(\frac{z-z_0}{2\lambda}\right)} >0. 
\end{gather*}
We consider $\mathbb{R}^{2}$ with the metric $\I=\psi^*\langle,\rangle$. Next, let us compute the Christoffel symbols. First, 
\begin{equation*}
\nabla_{\partial_y}\partial_y=\Gamma _{11}^{1}\partial_y+\Gamma
_{11}^{2}\partial_z.
\end{equation*}
Taking inner product with $\partial_y$, $\partial_z$ we obtain:
\begin{eqnarray*}
E\Gamma _{11}^{1}+F\Gamma _{11}^{2} &=&\I\left( \nabla_{\partial_y}\partial_y,\partial_y\right) =\frac{1}{2}\partial_y\left( \I\left( \partial_y,\partial_y\right) \right) =\frac{1}{2}E_{y}=0, \\
F\Gamma _{11}^{1}+G\Gamma _{11}^{2} &=&\I\left( \nabla_{\partial_y}\partial
y,\partial_z\right) =\partial_yF-\I\left( \partial_y,\nabla_{\partial_y}\partial_z\right) 
=-\I\left( \partial_y,\nabla_{\partial_z}\partial_y\right)\\ & = & -\frac{1}{2}
\partial_zE=0. 
\end{eqnarray*}
The solution to this system is $\Gamma _{11}^{1}=\Gamma _{11}^{2}=0$. Similarly, we compute the other Christoffel symbols:
\begin{gather*}
\Gamma _{12}^{1}=\Gamma _{12}^{2}=0, \quad 
\Gamma_{22}^1=\dfrac{-1}{4\lambda^2}, \quad \Gamma_{22}^2=
\dfrac{1}{2\lambda}\tanh\left(\dfrac{z-z_0}{2\lambda}\right).
\end{gather*}
Now the equations of a geodesic $\alpha \left( t\right) =\left( y(t),z(t) \right)$ are (cf. \cite[p. 67]{ON})
\begin{eqnarray*}
0 &=&y''+\Gamma _{11}^{1}(y')^{2}+2\Gamma _{12}^{1}y'z'+\Gamma _{22}^{1}(z') ^{2}
=y''-\dfrac{\left( z'\right)^{2}}{4\lambda ^{2}}, \\
0 &=&z''+\Gamma _{11}^{2}(y')^{2}+2\Gamma _{12}^{2}y'z'+\Gamma _{22}^{2}(z')^{2}
 \\ &=& z''+\dfrac{1}{2\lambda} \tanh \left( \dfrac{z-z_0}{2\lambda }\right) ( z') ^{2}.
\end{eqnarray*}
The general solution to the second ODE is
\[ z(t)=z_0+2\lambda\,\mathrm{asinh}\left(\dfrac{a_1}{2\lambda}t+a_2\right), \quad a_1,a_2\in\R,
\]
where $\mathrm{asinh}:\R\to\R$ is the globally defined inverse function of $\sinh$. Moreover,
\[ y''(t) = \frac{a_1^2}{(a_1t+2\lambda a_2)^2+4\lambda^2}.
\]
Integrating here, we obtain
\[ y'(t)=\frac{a_1}{2\lambda} \arctan\left(\frac{a_1 t+2a_2\lambda}{2\lambda}\right)+b_1.
\]
As $\vert\arctan(x)\vert <\pi/2$ for any $x\in\R$, then for a suitable constant $A>0$, $\vert y'(t)\vert \leq A$ for any $t\in\R$. Therefore, any inextensible solution $y$ is globally defined on the whole $\R$. 
\end{pf}
\begin{rmk} \normalfont \label{novan}
Recall that a properly immersed space-like hypersurface in Min\-kowski $n$-space whose normal vector satisfies the subaffine growth condition is complete (see \cite{BE}). Also,  if a properly immersed space-like hypersurface in Min\-kowski $n$-space has bounded principal curvatures, then it is complete (see  \cite{turca} and \cite{stacey}). 

Take an inextensible example of type \ref{caso1} This entire graph is properly embedded. The partial derivatives of the immersion are
\[  \psi_y = (-2\lambda^2,1,0), \ \psi_z=\left(2\lambda  \tanh\left(\frac{z-z_0}{2\lambda}\right), 0,1 \right). 
\]
The coefficients of the first fundamental form are
\[ E=4\lambda^2, \ F= -2\lambda \tanh\left(\frac{z-z_0}{2\lambda}\right), \ G=1. 
  \]
The normal vector is
\[ N=\left( \lambda \cosh\left(\frac{z-z_0}{2\lambda}\right), \frac{1}{2\lambda}\cosh\left(\frac{z-z_0}{2\lambda}\right), 
\sinh\left(\frac{z-z_0}{2\lambda}\right)\right).
\]
However, this vector does not satisfy the subaffine growth condition, because its coordinates behave as the exponential map at infinity. 
The coefficients of the second fundamental form are
\[ e= f=0, \ g=\frac{-1}{2\lambda}\mathrm{sech}\left(\frac{z-z_0}{2\lambda}\right). \]
From \eqref{shape}, a straightforward computation  gives  the principal curvatures $\lambda_1=0$ and
\[\lambda_2= \frac{-1}{2\lambda}\cosh\left(\frac{z-z_0}{2\lambda}\right).\] 
Clearly, this function is not bounded. $\Box$
\end{rmk}

\begin{corollary}\label{tipo-II-completeness}
\label{type2} All inextensible solutions of type  $III$ %\ref{caso3} 
in Theorem \ref{convolution-surfaces} are time-like, never entire, and incomplete (space-like, time-like, light-like).
\end{corollary}
\begin{pf} Given $b_{0}$, $z_{0}\in \mathbb{R}$, $\lambda >0$, we take $u(y,z)=2\lambda ^{2}y-4\lambda ^{2}\log_e \left\vert \cos	\left( \frac{z-z_0}{2\lambda }\right) \right\vert$ $+b_{0}$, and $\psi(y,z) =\left( u(y,z),y,z\right)$. Firstly, they cannot be entire because they can only be defined on horizontal strips of the form
\begin{equation*}
S(z_0,\lambda,k)=\left\{ (y,z)\in \mathbb{R}^2 :  
-\frac{\pi}{2}+k\pi <
\frac{z-z_0}{2\lambda} < \frac{\pi}{2}+k\pi  
\right\}, \ k\in\mathbb{Z}.
\end{equation*}
We recall the coefficients of the first fundamental form: 
\begin{equation*}
\left(\begin{matrix} E & F \\ F & G\end{matrix}\right)=
\left(\begin{matrix} -2u_y & -u_z \\ -u_z & 1\end{matrix}\right) 
=\left( \begin{matrix}
-4\lambda ^{2} & -2\lambda \tan \left( \frac{z-z_0}{2\lambda }\right) \\ 
-2\lambda \tan \left( \frac{z-z_0}{2\lambda }\right) & 1
\end{matrix}\right).
\end{equation*}
The causal character of $\psi ^{\ast }g$ is determined by  
\begin{equation*}
EG-F^2=-4\lambda ^{2}\left( 1+\tan ^{2}\left( \frac{z-z_0}{2\lambda}
\right) \right)<0,
\end{equation*}
namely   the surface $\psi$ is time-like. Take the metric $\I=\psi^*\langle,\rangle$.  Since  $\psi:(S(z_0,\lambda,k),\I)\to \LL$ is an isometric embedding, the map $\psi$ is an isometry onto its image. Then  we can work on $(S(z_0,\lambda,k),\I)$. This surface is simply connected, and $\I$ is a Lorentzian metric. By Corollary \ref{flat}, $(S(z_0,\lambda,k),\I)$ is flat, so that it is globally isometric to an open subset of the Minkowski plane, say $\Phi:(\Omega,g_o)\to (S(z_0,\lambda,k),\I)$, where $\Omega\subset \L^2$ and $g_o$ is the standard metric on $\L^2$. 

With this information, given $k\in\mathbb{Z}$, we consider the following curve $\alpha:(-\pi\lambda,\pi\lambda)\rightarrow S(z_0,\lambda,k)\subset \mathbb{R}^2$,  
$\alpha(t)=(y_0,t+z_0+2k\lambda\pi).$ This curve is inextensible and divergent. Simple computations show  $\vert\alpha'(t)\vert^2=1$, so its total length is $L(\alpha)=2\pi\lambda$. But now, since $\Phi$ is an isometry, the curve $\beta=\Phi^{-1}\circ\alpha$ is inextensible, space-like, unit, divergent, with total length $2\pi\lambda$. As $(\Omega,g_o)$ is an open subset of $\L^2$, then $\Omega$ is not the whole $\L^2$. This readily shows that the surface is not complete (in any sense, space-like, time-like, light-like).
\end{pf}
\begin{corollary}\label{tipo-III}
Any inextensible solution of type $IV$ %\ref{caso4} 
in Theorem \ref{convolution-surfaces} is
time-like, not entire, and not complete (space-like, time-like, light-like).
\end{corollary}

\begin{pf} Given $a_{0},z_{0}\in \mathbb{R}$, $\lambda>0$, let us consider 
$u(y,z)$ as in case $IV$ %\ref{caso4} 
and $\psi (y,z) =\left( u(y,z),y,z\right)$. Clearly,  they can only be defined on horizontal strips of the form 
	\begin{equation*}
		S^{+}(z_{0})=\left\{ (y,z)\in \mathbb{R}^{2}:z>z_{0}\right\} ,\
		S^{-}(z_{0})=\left\{ (y,z)\in \mathbb{R}^{2}:z<z_{0}\right\} .
	\end{equation*}
We compute  
\begin{gather*}
u_y=-2\lambda^2, \quad 
u_z=2\lambda \coth \left( \frac{z-z_{0}}{2\lambda }\right), \\
E=-2u_y, \ F=-u_z, \ G=1, \ EG-F^2 =\frac{-4\lambda ^{2}}{\sinh ^{2}\left( \frac{z-z_{0}}{2\lambda }\right) }<0. 
\end{gather*}
Next, the divergent curve $\alpha:(z_0,z_0+1)\rightarrow S^+(z_0)$, $\alpha(s)=(y_0,s)$,  satisfies that $\alpha'(s)=(0,1)$, it is inextensible at $z_0$, and $\vert\alpha'(s)\vert^2=1$. In particular, its total length is finite. Similarly, $\alpha:(z_0-1,z_0)\rightarrow S^-(z_0)$, $\alpha(s)=(y_0,s)$ is an inextensible, divergent curve with finite total length. By repeating the argument in Corollary \ref{tipo-II-completeness}, we obtain that this surface cannot be complete. 
\end{pf}

\begin{corollary}\label{type1} Any solution $\psi$ of  $II$ %type \ref{caso2} 
in Theorem is space-like when $a_1>0$, and time-like when $a_1<0$. Moreover, any inextensible solution $\psi$ is entire, but not complete (space-like, time-like, light-like). 
\end{corollary}

\begin{pf}
We consider $\psi (y,z) =\left( u(y,z),y,z\right)$, with $u(y,z)  =a_{1}e^{-2y}-\dfrac{b_{1}^{2}}{2}y+b_{1}z+b_{0}$, for $a_{1},b_{1},b_{0}\in \mathbb{R}$. 	First, it is very clear the $u$ can always be extended to the whole $u:\mathbb{R}^2\rightarrow\mathbb{R}$. We compute the coefficients of the first fundamental form: 
\begin{equation*}
\left(\begin{matrix} E & F \\ F & G\end{matrix}\right)=\left( 
		\begin{matrix} 	-2u_y & -u_z \\ -u_z & 1 \end{matrix}
		\right) =\left( 
		\begin{matrix}
			b_{1}^{2}+4a_{1}e^{-2y} & -b_{1} \\ 
			-b_{1} & 1
		\end{matrix}
		\right).
\end{equation*}
Clearly, the surface is space-like for $a_{1}>0$, namely  the metric $\I$ is Riemannian. Vice versa, the surface is time-like when $a_1<0$. Let us consider $\mathbb{R}^{2}$ with the metric $\I=\psi^*\langle,\rangle$. We compute the Christoffel  symbols:
\[ \Gamma_{11}^1=-1, \quad \Gamma_{11}^2=-b_1, \quad \Gamma_{12}^1=\Gamma_{12}^2=\Gamma_{22}^1=\Gamma_{22}^2=0.
\]
The equations of a geodesic $\alpha \left( t\right) =\left( y\left(
t\right),z\left( t\right) \right)$ are
\begin{eqnarray*}
0 &=&y''+\Gamma _{11}^{1}(y')^{2}+2\Gamma _{12}^{1}y'z'+\Gamma _{22}^{1}(z') ^{2}
=y''-(y')^2, \\
0 &=&z''+\Gamma _{11}^{2}(y')^{2}+2\Gamma _{12}^{2}y'z'+\Gamma _{22}^{2}(z')^{2}
=z''-b_1(y') ^{2}.
\end{eqnarray*}
It is simple to check that the divergent curve 
\[\alpha:[0,1)\to \R^2,\quad \alpha(t)=\left(\log_e(1-t),b_1\log_e(1-t)\right)\] is a geodesic. By a simple computation, we obtain  $\I_{\alpha(t)}\big(\alpha'(t),\alpha'(t)\big) = 4a_1. $ 
This implies that the length of $\alpha$ is $\mathrm{length}(\alpha)=\int_o^1 \vert\alpha'(t)\vert dt = 
2\sqrt{\vert a_1\vert}$. In either case ($a_1>0$ or $a_1<0$), there exists a divergent geodesic with finite total length. This means that this geodesic is not complete. By repeating the argument of Corollary \ref{tipo-II-completeness}, this surface is not complete (space-like, time-like, light-like).
\end{pf}
\begin{example}\normalfont  Plenty of mean curvature flows in the Minkowski  2-plane are computed in \cite{L2}. In particular, there are essentially 3 translating solitons, travelling in space-like, time-like and light-like directions. These curves can be seen as the corresponding \textit{Grim Reaper} curves in $\mathbb{L}^2$. Those in a light-like direction are written in our coordinates as $x=e^{-2y}/2$.  If we consider the projection $\pi:\mathbb{L}^3\to \mathbb{L}^2$, $\pi(x,y,z)=(x,y)$, the preimage by $\pi$ of each curve provides a translating soliton in a space-like, time-like and light-like direction, respectively. We will regard them as the \textit{Grim Reaper} surfaces in $\mathbb{L}^3$. If we take $a_1=2$, $b_1=0$ in case $II$ 
in Theorem \ref{convolution-surfaces}, we obtain the corresponding \textit{Grim Reaper} surface in a light-like direction, but written our way. $\Box$
\end{example}
We summarize these corollaries in the following table:
\begin{table}[h]
\begin{center}
\begin{tabular}{cccc}
Type & Entire & Causal Character & Completeness \\ \hline
I & yes & space-like & yes \\ \hline 
II & yes & space-like if $a_1>0$  & no \\
 & yes & time-like if $a_1<0$ & no \\ \hline  
III & no & time-like & no  \\ \hline 
IV & no & time-like & no \\ \hline 
\end{tabular}
\end{center}
\caption{Causal character and completeness.}
\end{table}

\section{The Group of Isometries Whose Axis Is Light-like}\label{eje}

We use the following subgroup of direct, time-orientation preserving isometries
\begin{equation*}
A_{3}=\left\{ \xi _{t}= 
\begin{pmatrix}
1 &  0 & 0 \\ \frac{1}{2}t^{2} & 1 & t \\ t & 0 & 1
\end{pmatrix} : t \in \mathbb{R}\right\} .
\end{equation*}
The action is given by $(x,y,z)\in\LL$, $\xi_{t}\cdot (x,y,z)=(x,y,z)\xi_{t}$, with the usual matrix multiplication. We will need the following regions in $\mathbb{L}^{3}\backslash \left\langle 
\vec{x}\right\rangle$: 
\begin{gather*}
\mathcal{S}^{+}=\left\{ \left( x,y,z\right) \in \mathbb{L}^{3}:y>0\right\},  \ 
\mathcal{S}^{-}=\left\{ \left( x,y,z\right) \in \mathbb{L}^{3}:y<0\right\},\\
\mathcal{S}=\left\{ \left( x,y,z\right) \in \mathbb{L}^{3}:z=0\right\},
\end{gather*}
and inside them, the following open half planes:
\begin{gather*}
\widetilde{\mathcal{S}}^{+} =\mathcal{S}^{+}\cap \mathcal{S}=
\left\{(x,y,0) \in \mathbb{L}^{3}:y>0\right\}, \\
\widetilde{\mathcal{S}}^{-} =\mathcal{S}^{-}\cap \mathcal{S}=\left\{
(x,y,0) \in \mathbb{L}^{3}:y<0\right\}.
\end{gather*}
We recall the following result from \cite{BCO}.
\begin{thm}\label{fundregions}
Let M be a connected surface and $\Phi :M\rightarrow \mathbb{L}^{3}$ a
non-degenerate immersion. Then  $(M,\Phi ^{\ast }(g))$ is $A_{3}$-invariant
if and only if one of the following statements holds:
\begin{enumerate}
\item If $(M,\Phi ^{\ast }(g))$ is Riemannian, there exists a regular space-like
curve $\alpha $, immersed in either $\widetilde{\mathcal{S}}^{+}$ or $\widetilde{\mathcal{S}}^{-}$, such that $\Phi \left( M\right) =\left\{ \xi
_{t}\left( trace\left( \alpha \right) \right) :t\in \mathbb{R}\right\} $.

\item If $(M,\Phi ^{\ast }(g))$ is Lorentzian, there exists a regular time-like
curve $\alpha $, immersed in either $\widetilde{\mathcal{S}}^{+}$ or $\widetilde{\mathcal{S}}^{-}$, such that $\Phi \left( M\right) =\left\{ \xi
_{t}\left( trace\left( \alpha \right) \right) :t\in \mathbb{R}\right\} $. 
\end{enumerate}
\end{thm}

Now we take a regular curve in $\mathcal{S}$, $\alpha :I\subseteq \mathbb{R}\rightarrow \mathcal{S}\subset \mathbb{L}^{3}$, $\alpha (s)=(x(s),y(s),0)$, and construct  the $A_{3}$ invariant surface $\psi (s,t)$ as:
\begin{equation*}
\psi:I\times\R\to \LL, \quad \psi (s,t) 
=\left( x(s) +\frac{t^{2}}{2}y(s),y(s),ty(s)\right).
\end{equation*}
We compute the partial derivatives of $\psi(s,t)$, 
\begin{equation*}
\psi _{s}=\left( x'(s) +\frac{t^{2}}{2}y'(s),y'(s),ty'(s)\right), \quad 
\psi _{t}=\left( ty(s),0,y(s)\right).
\end{equation*}
The matricial expression of the first fundamental form $\I=\psi^{*}\langle,\rangle$ is 
\begin{equation*}
\left(\begin{matrix} E & F \\ F & G \end{matrix} \right) = 
\left(\begin{matrix}
-2x'(s) y'(s) & 0 \\ 
0 & y^{2}(s)
\end{matrix}\right). 
\end{equation*}
Since we assume that the induced metric is not degenerate, we obtain 
 $x'(s) y'(s) \neq 0$. Then  we define 
\begin{equation*}
\varepsilon =\mathrm{sign}\left( x'(s) y'(s) \right) =\pm 1,\quad 
W=\sqrt{\dfrac{2\varepsilon x'(s) }{y'(s) }}>0,
\end{equation*}
and we construct the unit normal vector 
\begin{equation*}
N=\frac{1}{W}\left( \frac{t^{2}}{2}-\frac{x'(s) }{y'(s) },1,t\right).
\end{equation*}
Note that $\langle N,N\rangle=\varepsilon$. The second partial derivatives of $\psi$ are
\begin{align*}
\psi_{ss}=&\left( x''(s) +\frac{t^{2}}{2}y''(s),y''(s),t y''(s)\right), \\
\psi_{st}=& \left( t y'(s),0,y'(s)\right), \quad 
\psi_{tt}=(y(s),0,0).
\end{align*}
The coefficients of the second fundamental form are
\[ e=\langle N,\psi_{ss}\rangle = \frac{x' y''-x'' y'}{W y'},\quad 
f =\langle N,\psi_{st}\rangle= 0, \quad g=\langle N,\psi_{tt}\rangle=\frac{-y}{W}. 
\]
Therefore 
\[ H = -\frac{y\, x'\, y''+2 x'\, (y')^2 -y\, y'\, x''}{2 W\, y\, x' \,(y')^2}.
\]
Next, by Definition \ref{Tlight-like}, with $\vec{K}=\vec{x}$ the
chosen parallel light-like vector field, we compute:
\begin{align*}
\langle \vec{x}^{\perp },N\rangle = &
\langle \vec{H},N\rangle =\langle \varepsilon HN,N\rangle =H, \\
\langle \vec{x}^{\perp },N\rangle = &
\Big\langle (1,0,0),
\frac{1}{W}\left( \frac{t^2}{2}-\frac{x'(s)}{y'(s)},1,t\right) 
\Big\rangle =-\frac{1}{W}. 
\end{align*}
As a result, we obtain the following ODE:
\begin{equation}
y\,x'\,y''+2x'\,(y')^2-y\,y'\,x''=2y\,x'(y')^2.
\end{equation}

We want to express the profile curve $\alpha$ as a graph. Firstly, let us assume $y(s) =s$ and we examine $x(s)$.  Then  our equation transforms into
\[ 2x'(s)-s\,x''(s)=2s\,x'(s).
\]
A standard computation shows the general solution to this ODE, 
\begin{equation*}
x(s) =a_{0}\left( 2s^{2}+2s+1\right) e^{-2s}+a_{1},\ 
a_{1},a_{0}\in \mathbb{R}.
\end{equation*}
However, we have to discard the case $a_0=0$ because $x'(s)=0$ for any $s$. Coming back, we see that $x'(s)=-4a_0s^2e^{-2s}$, so that
\[\varepsilon = \mathrm{sign}(x'(s)y'(s)) = \mathrm{sign}(-a_0).
\]
Therefore, the normal $N$ is time-like when $\varepsilon=-1$, that is to say, when $a_0>0$. And it is space-like when $a_0<0$. In other words, $\psi$ is space-like when $a_0>0$, and time-like when $a_0<0$.

Now, we  assume $x(s) =s$ and let us examine $y(s)$. Our ODE transforms into
\[ y\,y''+2(y')^2=2y\,(y')^2. 
\]
Clearly,  $y(s)=y_0\in\R$ is a solution to this equation,  but then we get $y'(s)=0$, and we supposed that $x^{\prime}(s) y'(s) \neq 0$. If we arrange the above equality, we get  
\begin{equation*}
\frac{y''(s) }{y'(s) }
=2\left( 1-\frac{1}{y(s) }\right) y'(s).
\end{equation*}
By integrating both sides, we obtain $y^{2}(s) e^{-2y(s) }y'(s)=b_{0}$, $b_{0}\in \mathbb{R}$, $b_{0}\neq 0.$
We define the function 
\begin{equation*}
\phi :\mathbb{R\rightarrow R},\ \phi(r) =\frac{-1}{4}\left(2r^{2}+2r+1\right) e^{-2r}.
\end{equation*}
Note that $\phi'(r) =r^{2}e^{-2r}\leq  0$. Moreover, $\phi'(r) =0$ if, and only if, $r=0$. Therefore, $\phi $ is injective. Next, 
\[
\lim_{r\to +\infty}\phi(r)=0, \quad \lim_{r\to -\infty}\phi(r)=-\infty,
\quad \phi(0)=\frac{-1}{4}.\]
This shows $\phi:\R\to(-\infty,0)$.  To have a well-defined surface, the profile curve cannot get out of $\widetilde{\mathcal{S}}^{+}$ or $\widetilde{\mathcal{S}}^{-}$, so that we need to exclude $-1/4$ from the interval $J$. Therefore, the solutions are 
\[
y:J\to \R, \ y(s) =\phi ^{-1}\left( b_{0}s+b_{1}\right),
\]
where $J\subseteq(-\infty, -b_1/b_0)\backslash\{-1/4\}$, if $b_0>0$, or $J=(-b_1/b_0,+\infty)\backslash\{-1/4\}$, if $b_0<0$. Now, 
\[\varepsilon =\mathrm{sign}(x'(s)y'(s))=\mathrm{sign}\left(
b_0\dfrac{e^{2y(s)}}{ y(s)^2 }\right)=\mathrm{sign}(b_0).
\]

\begin{thm}\label{parabolic}
Let $\psi:I\times\R\to\LL$, $\psi=\psi (s,t)$, be an $A_{3}$ invariant surface, such that it is a translating soliton on the parallel light-like direction $\vec{K}=(1,0,0)$. Then the profile curve $\alpha :I\subset \mathbb{R}\longrightarrow \mathbb{L}^{3}$, $\alpha (s) =\left( x(s),y(s),0\right)$ of $\psi(s,t)$, is one of the following:

\begin{enumerate}
\item For $y(s) =s(\neq 0)$, given $a_{0},a_{1}\in \mathbb{R}$, $a_{0}\neq 0$, $x(s) =a_{0}\left( 2s^{2}+2s+1\right) e^{-2s}+a_{1}$. In addition, $\psi$ is space-like iff $a_{0}>0$, and time-like iff $a_{0}<0$. 

\item For $x(s) =s$, given $b_{2},b_{3}\in \mathbb{R}$, $b_{2}\neq 0$, and the diffeomorphism $\phi :\mathbb{R}\to (0,+\infty)$, $\phi(r)=\left( 2r^{2}-2r+1\right) e^{-2r}$, 
\begin{equation*}
y:J\rightarrow \mathbb{R},\ y(s) =\phi ^{-1}\left(-4 b_{2}s+b_{3}\right), 
\end{equation*}
where the interval $J$ is included in either $J\subseteq (-\infty, -b_1/b_0)\backslash\{1\}$, if $b_0>0$, or $J\subseteq (-b_1/b_0,+\infty)\backslash\{1\}$, if $b_0<0$. In addition, $\phi$ is space-like iff $b_0<0$, and $\psi$ is time-like iff $b_0>0$.
\end{enumerate}
\end{thm}

\begin{rmk}\normalfont 
When the surface approaches the affine plane $\mathcal{S}$, there are singularities. Indeed, case 1), for each $t\in\R$, $\lim_{s\to 0}\psi(s,t)=\big(a_0+a_1,0,0\big)$. \end{rmk}

\section*{Acknowledgements}

The second author is partially supported by the Spanish Ministry of Economy and Competitiveness, and European Region Development Fund, project MTM2016-78807-C2-1-P, and by the Junta de Andaluc\'{\i}a grant A-FQM-494-UGR18, with FEDER funds. The authors would like to thank the referees for their suggestions and questions, which improved the paper.

\end{document}